\documentclass[12pt]{amsart}

\usepackage{amssymb,epsf}

\newtheorem{thm}{Theorem}[section] 
\newtheorem{defn}[thm]{Definition}

\newtheorem{lem}[thm]{Lemma}
\newtheorem{exa}[thm]{Example}
\newtheorem{prop}[thm]{Proposition}

\newtheorem{observ}[thm]{Observation}

\newcommand{\R}{\mathbb R}

\begin{document}

\title[Orchard crossing number of complete bipartite graphs]{On the Orchard crossing number of complete bipartite graphs}

\author{Elie Feder and David Garber}

\address{Kingsborough Community College of CUNY, Department of Mathematics and Computer Science,
2001 Oriental Blvd., Brooklyn, NY 11235, USA}
\email{eliefeder@gmail.com, efeder@kbcc.cuny.edu}

\address{Department of Applied Mathematics, Faculty of Science, Holon Institute of Technology, Golomb 52,
PO Box 305, Holon 58102, Israel, and (Sabbatical:) Department of Mathematics, Jerusalem College of Technology, Jerusalem, Israel}
\email{garber@hit.ac.il}

\date{\today}

\begin{abstract}
We compute the Orchard crossing number, which is defined in
a similar way to the rectilinear crossing number, for the complete bipartite graphs $K_{n,n}$.
\end{abstract}

\maketitle

\section{Introduction}

Let $G$ be an abstract graph. Motivated by the Orchard relation,
introduced in \cite{bacher,BaGa}, we have defined the {\it Orchard
crossing number} of $G$ \cite{FG}, in a similar way to the well-known {\it
rectilinear crossing number} of an abstract graph $G$ (denoted by
$\overline{\rm cr}(G)$, see \cite{aak2,PT}). A general reference for crossing numbers
can be \cite{F}.

The Orchard crossing number is interesting for several reasons.
First, it is based on the Orchard relation which is an equivalence
relation on the vertices of a graph, with at most two equivalence
classes (see \cite{bacher}). Moreover, since the Orchard relation
can be defined for higher dimensions too (see \cite{bacher}),
hence the Orchard crossing number
may be also generalized to higher dimensions.

Second, a variant of this crossing
number is tightly connected to the well-known rectilinear crossing number (see
Proposition \ref{max_complete} below).

Third, one can find real problems which the Orchard crossing number
can represent.  For example, design a network of computers which should be
constructed in a manner which allows possible extensions
of the network in the future. Since we want to avoid (even future) crossings
of the cables which are connecting between the computers, we need to count not
only the present crossings, but also the separators
(which might come to cross in the future).

\medskip

In this paper, we compute the Orchard crossing number for the complete
bipartite graphs $K_{n,n}$:
\begin{thm}\label{main_thm}
$${\rm OCN}(K_{n,n})= 4n{n \choose 3}$$
This value is attained where all the $2n$ points are in a convex position and alternate in color (see Figure \ref{K44} for an example for $n=4$).

\begin{figure}[h]
\epsfysize=4cm
\centerline{\epsfbox{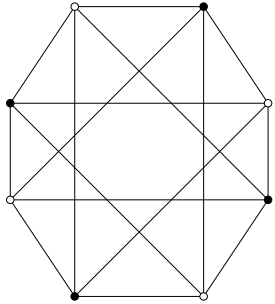}}
\caption{The optimal rectilinear drawing for ${\rm OCN}(K_{4,4})$} \label{K44}
\end{figure}

\end{thm}

The ideas of the proof are quite similar to those of \cite{AFH}, where the maximal value of the maximum rectilinear crossing number has been computed for some families of graphs, but still are not straightforward from them. This again shows the tight connection between the Orchard crossing number and the rectilinear crossing number.

The proof is based on two parts: in the first part, we show that the optimal drawing of $K_{n,n}$, presented in Figure \ref{K44} for $n=4$, has $4n{n \choose 3}$ Orchard crossings, and hence ${\rm OCN}(K_{n,n})\leq 4n{n \choose 3}$ (see Section \ref{upper}). In the second part, we show that any drawing of $K_{n,n}$ has at least $4n{n \choose 3}$ Orchard crossings, so ${\rm OCN}(K_{n,n})\geq 4n{n \choose 3}$ (see Section \ref{lower}).

\medskip

The paper is organized as follows.
In Section \ref{defs}, we present the Orchard relation, define
the Orchard crossing number, and give some examples.
Sections \ref{upper} and \ref{lower} contain the two parts of the proof of
Theorem \ref{main_thm}.

\section{The Orchard crossing numbers}\label{defs}

We start with some notations.
A finite set $\mathcal P = \{P_1,\dots, P_n\}$ of $n$ points
in the plane $\R^2$ is a {\it generic configuration}
if no three points of $\mathcal P$ are collinear.

A line $L \subset \R^2$ {\it separates} two points $P,Q \in
(\R^2\setminus L)$ if $P$ and $Q$ are in different connected
components of $\R^2 \setminus L$. Given a generic configuration
$\mathcal P$, denote by $n(P,Q)$ the number of lines defined by
pairs of points in $\mathcal P \setminus \{P,Q\}$, which separate
$P$ and $Q$.

\medskip

For defining the {\it Orchard crossing number} of an abstract graph $G$, we need some more notions.

\begin{defn}[Rectilinear drawing of an abstract graph $G$]
Let $G=(V,E)$ be an abstract graph, where $V$ is its set of vertices and $E$ is its set of edges.
{\em A rectilinear drawing of the abstract graph $G$}, denoted by $R(G)$, is a generic configuration of points $V'$
in the affine plane, in bijection with $V$. An edge $(s,t)\in E$ is
represented by the straight segment $[s',t']$ in $\R^2$.
\end{defn}

Then, we associate a {\it crossing number} to such a drawing:
\begin{defn}
Let $R(G)$ be a rectilinear drawing of the abstract graph $G=(V,E)$.
The {\em crossing number} of $R(G)$, denoted by $n(R(G))$, is:
$$n(R(G)) = \sum _{(s,t) \in E} n(s,t)$$
\end{defn}

Note that the sum is taken only over the edges of the graph, whence $n(s,t)$
counts in {\it all} the lines generated by pairs of points of the configuration.

\medskip

Now, we can define the {\it Orchard crossing number} of an abstract
graph $G=(V,E)$:

\begin{defn}[Orchard crossing number]
Let $G=(V,E)$ be an abstract graph.
The {\em Orchard crossing number} of $G$, denoted by ${\rm OCN}(G)$, is
$${\rm OCN}(G) = \min _{R(G)} (n(R(G)))$$
\end{defn}

\medskip

A variant of the Orchard crossing number is the {\it maximal Orchard crossing number}:

\begin{defn}[Maximal Orchard crossing number]
Let $G=(V,E)$ be an abstract graph. The {\em maximal Orchard crossing
number} of $G$, denoted by ${\rm MOCN}(G)$, is
$${\rm MOCN}(G) = \max _{R(G)} (n(R(G)))$$
\end{defn}

This variant is extremely interesting due to the following result (see \cite[Proposition 2.7]{FG}):

\begin{prop}\label{max_complete}
The rectilinear drawing which yields the maximal Orchard crossing number for
complete graphs $K_n$ is the same as the rectilinear drawing which attains the
rectilinear crossing number of $K_n$.
\end{prop}

The importance of this result is that it might be possible that the
computation of the maximal Orchard crossing number will be easier than
the computation of the rectilinear crossing number.

\section{An upper bound for ${\rm OCN}(K_{n,n})$}\label{upper}
In this section, we show the easy part of Theorem \ref{main_thm} by proving that the mentioned value $4n{n \choose 3}$ is indeed attained by the
rectilinear drawing of $K_{n,n}$, where all the $2n$ points are in a convex position and alternate in color (see Figure \ref{K44} for $n=4$).

\begin{lem} \label{num_comp_bipar}
Assume that $R(K_{n,n})$ is the rectilinear drawing of $K_{n,n}$ which
realizes the $2n$ points of $K_{n,n}$ as the points of a regular $2n$-gon,
and the points change colors alternately. Then:
$$n(R(K_{n,n}))= 4n{n \choose 3}$$
\end{lem}

\begin{proof}
Since any quadruple of points is in a convex position,
we have to consider only four types of quadruples:
\begin{enumerate}

\item All the points of the quadruple have the same color.
This case contributes nothing to the number of crossings (see Figure
\ref{bipartite-nn}(a)).

\item The quadruple consists of two black points and two white points,
and the points change colors alternately. This case contributes
nothing to the number of crossings (see Figure
\ref{bipartite-nn}(b), where the solid lines are edges of the graph and the dashed lines
are lines generated by a pair of points of the same color which are not edges of the graph).

\item Three of the points are black and one point is white (or vice versa).
This case contributes $1$ to the number of crossings (see Figure
\ref{bipartite-nn}(c)).

\item The quadruple consists of two consecutive black points followed
by two consecutive white points. This case contributes $2$ to the
number of crossings (see Figure \ref{bipartite-nn}(d)).

\end{enumerate}

\begin{figure}[h]
\epsfysize=2.5cm \centerline{\epsfbox{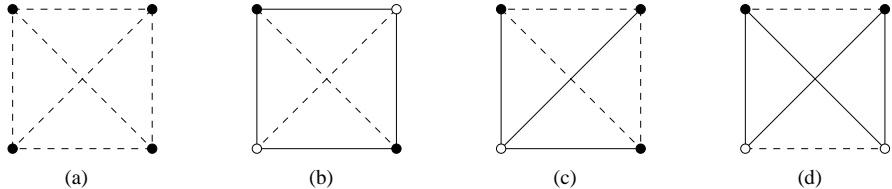}}
\caption{Four cases of quadruples of points in a convex position}
\label{bipartite-nn}
\end{figure}

Hence, we have to compute the number of quadruples of types (3) and
(4) respectively, and to multiply these numbers by their corresponding contributions
to the total number of crossings.

Assume that we have $n$ white points and $n$ black points.
For computing the number of quadruples of type (3), we have to choose
three black points out of $n$ black points, and then to choose one point
out of $n$ white points. There are $n{n \choose 3}$ possibilities to do so.
Since, we have to do the same with the opposite colors, we have
$2n{n \choose 3}$ quadruples of type (3),
which contribute $2n{n \choose 3}$ to the total number of crossings.

The count of quadruples of type (4) is a bit more complicated. We
count them in the following way: Choose an arbitrary point. Then
choose another point of the same color. Then, choose two points of
the opposite color to the left of the second point, but to the right of
the first point. One can easily see that the number of possibilities
is:
$$\frac{2n}{4} \sum_{k=2}^{n-1}k(k-1),$$
where the sum is induced from choosing the second point of the first
color. The division by $4$ comes from the fact that by this counting
argument, each quadruple is counted $4$ times, once for each point
of it.

We simplify this expression:
\begin{eqnarray*}
\frac{2n}{4} \sum_{k=2}^{n-1}k(k-1) & = & \frac{n}{2}
\sum_{k=1}^{n-1}(k^2-k) =\\
& = & \frac{n}{2}\left( \frac{(n-1)n(2n-1)}{6} -
\frac{n(n-1)}{2}\right) =\\
& =& \frac{n \cdot n (n-1)(n-2)}{6}=n{n \choose 3}
\end{eqnarray*}

Hence, we get that the number of quadruples of type (4) is $n{n
\choose 3}$. Since each quadruple of type (4) contributes two
crossings, this number should be multiplied by $2$ for getting the
total contribution of the quadruples of type (4). Hence, the
quadruples of type (4) contribute $2n{n \choose 3}$.

Summing up the two contributions yields the result.
\end{proof}

\section{A lower bound for ${\rm OCN}(K_{n,n})$}\label{lower}

In this section, we show the difficult part of Theorem \ref{main_thm} by proving that for any rectilinear drawing of $K_{n,n}$ there are at least
$4n{n \choose 3}$ Orchard crossings. This will show that $4n{n \choose 3}$ is indeed a lower bound for ${\rm OCN} (K_{n,n})$.

The idea of proof is counting separately the Orchard crossings induced by pairs of points with different colors and by pairs of points with the same color (see Sections \ref{first}
and \ref{second} respectively).

Together with the result of the previous section, this will prove Theorem \ref{main_thm} (see Section \ref{third}).

\medskip

We start with one notation. A {\em {\bf bw}-pair} is a pair of points consists of a black point and a white point.

\subsection{Orchard crossings induced by lines generated by pairs of points with different colors}\label{first}

Let $D$ be a configuration of $n$ white points and $n$ black points. For each
$k =1,2,\dots,n^2$, let $\ell_k$ be a line determined by a white point and a
black point. For each $\ell_k$, let $a_k$ be the number of {\bf bw}-pairs with one point in one
halfplane determined by $\ell_k$, and the other point in the other halfplane. Let $A= \sum\limits_{k=1}^{n^2} a_k$.

\begin{prop}\label{A}
$$A \geq 2n{n \choose 3}$$
\end{prop}

\begin{proof}
For each $\ell_k$, a black (resp. white) endvertex will be of {\em type $i$} if the
edge incident to $\ell_k$ divides the graph into two halfplanes, one
contains $i$ black (resp. white) vertices, and the other contains
$n-i-1$ black (resp. white) vertices.  By symmetry, we only have to consider
$0\leq i \leq \lfloor \frac{n-1}{2} \rfloor = N$.
Let $y_i$ be the number of endvertices of type $i$.  Thus, we have:
 $$y_0 + y_1+ \cdots +y_N = 2n^2,$$
since we have $n^2$ such edges and each one is counted twice
(for the white vertex and for the black vertex).

An edge which connects black and white vertices is of
{\em type $i,j$} if one halfplane determined by that edge has $i$
vertices of one color and $j$ vertices of the other color.  Let
$x_{i,j}$ be the number of edges of type $i,j$.
By symmetry, we can assume that $0 \leq i \leq j \leq N$
(we will justify this assumption later on). Note that $x_{i,j}=x_{j,i}$.

Thus, $y_i$ is related to $x_{i,j}$ by the following equation:
\begin{equation}\label{eqn_1}
y_i=2x_{i,i}+\displaystyle\sum_{j=0}^{i-1} x_{j,i} +\displaystyle\sum_{j=i+1}^{N} x_{i,j},
\end{equation}
due to the following argument: for being counted in $y_i$, an edge should have $i$ vertices
of one color in one of the halfplanes it determines. The only case which is counted twice is when $i=j$,
where it is counted for both colors.

Now, for an edge of type $i,j$, there are  $i(n-j-1)+j(n-i-1)$ {\bf bw}-pairs in opposite
halfplanes of that edge. Summing it over all edges of a drawing, we obtain
$$M=\displaystyle\sum_{i=0}^{N}\displaystyle\sum_{j=i}^{N} \left[ i(n-j-1)+j(n-i-1) \right] x_{i,j}$$
{\bf bw}-pairs in opposite halfplanes. We are looking for a drawing which minimizes $M$.

\medskip

We now justify our assumption that $0 \leq i \leq j \leq N $, where $i$ and $j$ are the number of
vertices of the two colors in the same halfplane of an edge.  Assume
that for a given type $i,j$ edge, the $i$ vertices of one color and
the $j$ vertices of the other color are in different halfplanes.
This yields $ij + (n-j-1)(n-i-1)$ {\bf bw}-pairs.  However,
$$i(n-j-1)+j(n-i-1) \leq ij + (n-j-1)(n-i-1)$$ for $0 \leq i \leq j \leq N $.
Hence our assumption reduces the number of {\bf bw}-pairs.

\medskip

In order to minimize $M$, we start by multiplying Equation (\ref{eqn_1}) by
$i(n-i-1)$, and subtracting it from $M$ for all values of $i$,
yielding:

\begin{eqnarray*}
& &M- \displaystyle\sum_{i=1}^{N} i(n-i-1)y_i =\\
&=&\displaystyle\sum_{i=0}^{N}\displaystyle\sum_{j=i}^{N} \left[ i(n-j-1)+j(n-i-1) \right] x_{i,j}\\
& &\qquad -\displaystyle\sum_{i=1}^{N} i(n-i-1)\left( 2x_{i,i}+\displaystyle\sum_{j=0}^{i-1} x_{j,i} +\displaystyle\sum_{j=i+1}^{N} x_{i,j}\right)=\\
& =&\displaystyle\sum_{i=0}^{N}\displaystyle\sum_{j=i+1}^{N} \left[ i(n-j-1)+j(n-i-1) \right] x_{i,j}\\
& & \qquad -\displaystyle\sum_{i=1}^{N} i(n-i-1)\left(\displaystyle\sum_{j=0}^{i-1} x_{j,i} +\displaystyle\sum_{j=i+1}^{N} x_{i,j}\right)=\\
& =&\displaystyle\sum_{i=0}^{N}\displaystyle\sum_{j=i+1}^{N} \left[ i(i-j)+j(n-i-1) \right] x_{i,j} \\
& & \qquad -\sum_{i=0}^{N-1}\sum_{j=i+1}^{N} j(n-j-1) x_{i,j}=\\
& =&\displaystyle\sum_{i=0}^{N}\displaystyle\sum_{j=i+1}^{N} \left[ i(i-j)+j(j-i) \right] x_{i,j} = \displaystyle\sum_{i=0}^{N-1} \displaystyle\sum_{j=i+1}^{N} (j-i)^2 x_{i,j}
\end{eqnarray*}

Hence, we have:
\begin{equation}\label{eqn_2}
M= \displaystyle\sum_{i=1}^{N} i(n-i-1)y_i +
\displaystyle\sum_{i=0}^{N-1} \displaystyle\sum_{j=i+1}^{N} (j-i)^2 x_{i,j}
\end{equation}

\medskip

For the next step, we introduce a new notation $p_{s,t}$. We first motivate it.
Note that every white (resp. black) vertex $v_i$ serves as an endvertex for $n$ edges of the graph. For each of these $n$ edges,  let $c_{i,j}$ ($1 \leq j \leq n$) be the number of black (resp. white) points in the halfplane with a smaller number of black (resp. white) points determined by this edge.
Thus, for each $v_i$, we have a sequence of $n$ numbers
$( c_{i,1}, \dots ,c_{i,n} )$ representing the types which $v_i$ is for the $n$ edges connected to $v_i$.

For example, for a white point $v$ on the convex hull of an alternating $2n$-gon (see Figure \ref{K44_example} for $n=4$, where only some of the edges are drawn), the two edges going to the adjacent black points (the edges $e_1,e_4$) are of type $0$. The next two edges (edges $e_2,e_3$) will be of type $1$, the next two edges (do not exist in this configuration) will be of type $2$, etc. Note that in this case we will have the same sequence of types for any point on the convex hull.

\begin{figure}[h]
\epsfysize=4cm
\centerline{\epsfbox{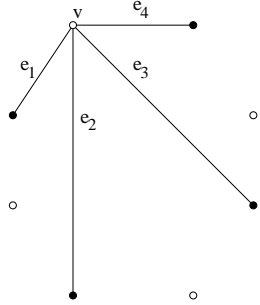}}
\caption{An example for computing $p_{s,t}$} \label{K44_example}
\end{figure}

Let $p_{s,t}$ be the number of white (resp. black) vertices having
$$s=\min\{ c_{i,1}, \dots ,c_{i,n} \},$$
and the index $t$ is the number of distinct sorted sequences generated by all the vertices.
For example, for the $2n$-gon configuration, we have $p_{0,1}= 2n$, since all $2n$ vertices have $0$ as the lowest term in the sequence $0,0,1,1,2,2,\dots$. Note that there is a unique sequence for all the vertices. Therefore, $p_{0,t}=0$ for $t \geq 2$, since there is no vertex whose lowest type is $0$ which has a different sequence. Moreover, $p_{s,t}=0$ for $s \geq 1$, since there is no sequence whose lowest term is greater than $0$.

Since $p_{s,t}$ counts number of vertices, and in total there are $2n$ vertices, we have:
$$\displaystyle\sum_{s=0}^{N} \displaystyle\sum_{t \geq 1} p_{s,t}=2n.$$


Denote by $z_{s,t,i}$ the number of appearances of $i$ in the $t$th sequence whose lowest term is $s$.
For example, assume that the first sequence is $0,1,2,1,0$. Then: $z_{0,1,0} =2$ because $0$ appears twice in the sequence. Similarly, $z_{0,1,1}=2$ since $1$ appears twice and $z_{0,1,2}=1$ since $2$ appears once.

It follows that:
\begin{equation}\label{eqn_3}
y_i=\displaystyle\sum_{t \geq 1} \displaystyle\sum_{s=0}^{i} z_{s,t,i}p_{s,t}
\end{equation}

Additionally, since every vertex has $n$ edges, for fixed $s$ and $t$ we have that
\begin{equation}\label{eqn_5}
\displaystyle\sum_{i=s}^{N} z_{s,t,i}=n.
\end{equation}

Since $\displaystyle\sum_{t \geq 1} \sum_{s=0}^{N} p_{s,t}=2n$, we obtain by Equation (\ref{eqn_3}):
\begin{equation}\label{eqn_6}
y_i= 4n +  \displaystyle\sum_{t \geq 1} \left[ \displaystyle\sum_{s=0}^{i} (z_{s,t,i}-2)p_{s,t} - 2\displaystyle\sum_{s=i+1}^{N}p_{s,t}\right].
\end{equation}
Respectively, for odd $n$, we have:
\begin{equation}
y_N= 2n +  \displaystyle\sum_{t \geq 1} \displaystyle\sum_{s=0}^{N} (z_{s,t,N}-1)p_{s,t} .
\end{equation}

We continue for even $n$. Using Equation (\ref{eqn_6}), we can rewrite the first
part of the expression for $M$ (Equation (\ref{eqn_2})) as follows:
\begin{eqnarray*}
\sum_{i=1}^{N} i(n-i-1)y_i & = & 4n \sum_{i=1}^{N} i(n-i-1)+\\
& & + \sum_{t \geq 1} \sum_{i=1}^{N} i(n-i-1)\left[\sum_{s=0}^{i}
(z_{s,t,i}-2)p_{s,t} - 2\sum_{s=i+1}^{N}p_{s,t}\right]
\end{eqnarray*}
Following a change in the indices of the sums, this can be rewritten
as:
$$4n \displaystyle\sum_{i=1}^{N} i(n-i-1) + \displaystyle\sum_{t\geq 1} \displaystyle\sum_{s=0}^{N}
p_{s,t}\left[\displaystyle\sum_{i=s}^{N} i(n-i-1)(z_{s,t,i}-2)-2\displaystyle\sum_{i=1}^{s-1} i(n-i-1)\right]$$
This can again be rewritten as:
\begin{eqnarray*}
& & 4n \sum_{i=1}^{N} i(n-i-1) + \left.\sum_{t\geq 1}\sum_{s=0}^{N} p_{s,t}\right[s(n-s-1)(z_{s,t,s}-2)+\\
& & \quad \left.+\sum_{i=s+1}^{N} i(n-i-1)(z_{s,t,i}-2)-2\sum_{i=1}^{s-1} i(n-i-1)\right]=\\
& = & 4n \sum_{i=1}^{N} i(n-i-1)+\sum_{t \geq 1} \sum_{s=0}^{N} p_{s,t}\left[ s(n-s-1)\sum_{i=s}^{N} (z_{s,t,i}-2)+\right.\\
& & +\left.\displaystyle\sum_{i=s+1}^{N}[i(n-i-1)-s(n-s-1)](z_{s,t,i}-2)-2\displaystyle\sum_{i=1}^{s-1} i(n-i-1)\right]
\end{eqnarray*}

Using Equation (\ref{eqn_5}), it follows that this is also equal to:
\begin{eqnarray*}
& & 4n \sum_{i=1}^{N} i(n-i-1) +\left. \sum_{t \geq 1} \sum_{s=0}^{N} p_{s,t}\right[C(s,n)+ \\
& & \quad + \left. \sum_{i=s+1}^{N}[i(n-i-1)-s(n-s-1)](z_{s,t,i}-2)\right],
\end{eqnarray*}
where
\begin{eqnarray*}
C(s,n)&=& s(n-s-1)(n - \displaystyle\sum_{i=s}^{N} 2) - 2\displaystyle\sum_{i=1}^{s-1} i(n-i-1)= \\
&=& s(n-s-1)(n-2(N-s+1))- 2\displaystyle\sum_{i=1}^{s-1} i(n-i-1)
\end{eqnarray*}

Now, we show that $C(s,n)$ is non-negative for all $0 \leq s \leq N$ and $n$. For the proof, we will use the
following observation:
\begin{observ}\label{obs1}
Let $a,b$ be such that: $0<a<b<\frac{n}{2}$. Then:
$$a(n-a) < b(n-b)$$
\end{observ}

\begin{lem}
For all $0 \leq s \leq N$ and $n$, $C(s,n) \geq 0$.
\end{lem}

\begin{proof}
By the definition of $N=\lfloor \frac{n-1}{2} \rfloor$, we have:
$$s(n-s-1)(n-2(N-s+1)) \geq s(n-s-1)(2s-1).$$

Moreover, by Observation \ref{obs1},
\begin{eqnarray*}
2\displaystyle\sum_{i=1}^{s-1} i(n-i-1) &\leq& 2\displaystyle\sum_{i=1}^{s-1}(s-1)(n-s) = \\
&=& 2(s-1)^2(n-s) \\
&<& s(n-s-1)(2s-2).
\end{eqnarray*}

Therefore,
$$C(s,n)>s(n-s-1)(2s-1)-s(n-s-1)(2s-2)=s(n-s-1) \geq 0.$$
\end{proof}

We now show that $z_{s,t,i} -2 \geq 0$. Since the term $z_{s,t,i}-1$ must be carried throughout this summation, the expression for odd $n$ is also minimized for $p_{s,t}=0$ for all $s \geq 1$, provided $z_{s,t,N} -1 \geq 0$.

\begin{lem}\label{lem_part1}
$z_{s,t,i} \geq 2$ for all $s,t,i$ such that $i>s$ and $z_{s,t,N} \geq 1$ for odd $n$.
\end{lem}

\begin{proof}
Without loss of generality, consider a given black vertex. We start
by proving that there is at least one endvertex of type $\frac{n-1}{2}$
for odd $n$, and at least two endvertices of type $\frac{n-2}{2}$
for even $n$.  This statement can be proved by induction on $n$.
This statement is obvious for $n=2$ and $n=3$, so we start
with the inductive step.  Also, note that in traversing the $n$
edges incident to the given vertex in a clockwise or
counterclockwise manner, in moving from edge to edge, the number
of white vertices in the clockwise following halfplane may be changed by
at most $1$.  This fact will be used numerous times throughout the
proof. We illustrate it by the following example:

\begin{exa}\label{example1}
Given a rectilinear drawing of $K_{4,4}$ at Figure \ref{config_example} (where only some of edges are drawn). If we are traversing clockwise the $4$ edges starting frow the lowest black point, we have that the type of the edge $e_1$ is $0$, since there is no white point in the left halfplane defined by this edge. Next, the type of $e_2$ is $1$, since there is only one white point in the left halfplane defined by this edge. The type of $e_3$ is again $1$, since there is only one white point in the {\em right} halfplane defined by this edge. Finally, the type of $e_4$ is $0$, since there is no white point in the right halfplane defined by this edge.

Hence, we have that while moving from edge to edge, the number
of white vertices in the clockwise following halfplane may be changed by
at most $1$. Consequently, the types of the corresponding edges may be changed by
at most $1$.

\begin{figure}[h]
\epsfysize=3cm
\centerline{\epsfbox{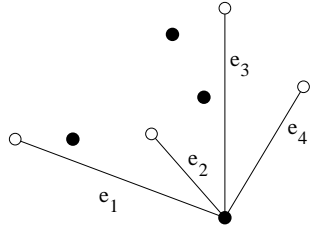}}
\caption{An example} \label{config_example}
\end{figure}

\end{exa}

\medskip

\noindent
\textbf{Case I:} Passing from odd $n$ to $n+1$. \\
Consider the edge for which this endvertex is of type
$\frac{n-1}{2}$ in the configuration of $n$ white vertices and $n$
black vertices. When the $(n+1)$st pair of vertices is added, this original
endvertex will be the first endvertex of type $\frac{(n+1)-2}{2}$.
If the $(n+1)$st white vertex is added in this edge's clockwise
following halfplane, then an immediately following edge or edge
extension has an endvertex of type $\frac{(n+1)-2}{2}$.  Thus, either
this edge or the edge corresponding to this extension will have the
second endvertex of type $\frac{n-1}{2}$.

\medskip

\noindent
\textbf{Case II:} Passing from even $n$ to $n+1$. \\
Consider an edge with an endvertex of type $\frac{n-2}{2}$ which has
$\frac{n}{2}$ white vertices in one of its halfplanes and
$\frac{n-2}{2}$ in the other.  If the $(n+1)$st white vertex is added
in the halfplane with $\frac{n-2}{2}$ vertices, then the considered
edge is now of type $\frac{(n+1)-1}{2}$.  If the $(n+1)$st vertex is added
in the halfplane with $\frac{n}{2}$ white vertices, then there are
$\frac{(n+1)+1}{2}$ white vertices in this halfplane and
$\frac{(n+1)-3}{2}$ white vertices in the clockwise following
halfplane of this edge's extension.  Since the number of white
vertices in the clockwise following halfplane can be changed by at most
$1$ when moving from edge to edge (edge ray and edge
extension), we find that traversing the graph from the edge with
$\frac{(n+1)+1}{2}$ white vertices in the clockwise following halfplane to
the extension with $\frac{(n+1)-3}{2}$, there must occur an edge or
extension with $\frac{(n+1)-1}{2}$ white vertices in the clockwise
following halfplane.  Thus, this edge or the edge corresponding to
the extension has an endvertex of type $\frac{(n+1)-1}{2}$. This completes
the proof for the maximal values.

\medskip

Using this result and the fact that in moving from edge to
adjacent edge, the number of white vertices in the clockwise
following halfplane may be changed by at most $1$, we can prove that
there are two endvertices of each type from the type $s+1$ to
the maximal type $N$.

\medskip

We split the proof according to the parity of $n$.

\begin{itemize}
\item For odd $n$, we have one endvertex of maximal
type $N=\frac{n-1}{2}$.  Traversing the $n$ edges starting and
ending with the edge with an endvertex of type $N$, from edge to
edge we must go down to an edge or an extension with $s$
vertices in the clockwise following halfplane, and then back up to
an edge with $N$.  Thus, we find there are at least two edges or
extensions with endvertices of each type from $s+1$ to $N$.

\item For even $n$, we have two edges with endvertices of maximal type
$N=\frac{n-2}{2}$.  Traversing the $n$ edges from one of the edges of type
$N$ to the other must go down to an edge or an extension with $s$
edges in the clockwise following halfplane and back up to an edge with
$N$.  Thus again, there are at least two edges or extensions with
endvertices of each type from $s+1$ to $N$.
\end{itemize}

\medskip

Hence, it follows that $z_{s,t,i} \geq 2$ for all $s,t,i$ such that $i>s$, and $z_{s,t,N} \geq 1$ for odd $n$ as needed.
\end{proof}

Additionally, by Observation \ref{obs1}, $\displaystyle\sum_{i=s+1}^{N}(i(n-i-1)-s(n-s-1)) \geq 0 $ for $i,s \leq
N=\frac{n-1}{2}$.

Going back to the final expression for Equation (\ref{eqn_2}), we have:

\begin{eqnarray*}
M&=&4n \sum_{i=1}^{N} i(n-i-1) + \left.\sum_{t \geq 1} \sum_{s=0}^{N} p_{s,t} \right[C(s,n)+\\
 & &+ \left.\sum_{i=s+1}^{N}(i(n-i-1)-s(n-s-1))(z_{s,t,i}-2)\right] +\\
 & &+ \sum_{i=0}^{N-1} \sum_{j=i+1}^{N} (j-i)^2 x_{i,j}
\end{eqnarray*}

Since $C(s,n)$, $(z_{s,t,i} -2)$ and $(j-i)^2$ are non-negative,
we find that this expression is minimized when
$p_{s,t}=0$ for $s \geq 1$ and $x_{i,j}=0$ for $i<j$.  For $s=0$, the expression is minimized when $z_{s,t,i}=2$ for all $i$. Evaluating
the sum for these conditions, we have:
$$M=4n \displaystyle\sum_{i=1}^{N} i(n-i-1) = 2n{n \choose 3}$$
{\bf bw}-pairs in opposite halfplanes determined by lines connecting two vertices
of opposite colors.

\medskip

By similar arguments, in the case where $n$ is odd, in the drawing for
which $M$ is minimized, we have the following expression:
$$M=4n \sum_{i=1}^{N-1} i(n-i-1) + 2n(N(n-N-1))= 2n{n \choose 3}.$$

\end{proof}

\subsection{Orchard crossings induced by lines generated by pairs of points of the same color}\label{second}

For each $k = 1,2,\dots, {n \choose 2}$, let $\ell_k$ be a line determined by two
points of the same color.
Without loss of generality, we assume that the two points are white.
The computation for a pair of black points is the same. Note that there is no edge in
$K_{n,n}$ based on this line, but still this line is counted within the separating lines.
Let $b_k$ be the number of {\bf bw}-pairs with one point
in one halfplane determined by $\ell_k$ and the other point in the other halfplane.
Let $B = \sum\limits_{k=1}^{n \choose 2} b_k$

\begin{prop}\label{B}
$$B\geq n{n \choose 3}$$
\end{prop}

\begin{proof}
According to our assumption, each $\ell_k$ is
determined by two white vertices. For each $\ell_k$, let an endvertex
be of {\em type $i$} if the line $\ell_k$ divides the graph into two
halfplanes, one containing $i$ white vertices, and the other
containing $n-i-2$ white vertices.  By symmetry, we have to consider only $0
\leq i \leq \lfloor \frac{n-2}{2} \rfloor = N $.

Let $y_i$ be the number of endvertices of type $i$.   Thus, we have
$$y_0 + y_1+ \cdots +y_N = 2{n \choose 2},$$
since we have ${n \choose 2}$ pairs of white vertices and each pair is counted twice for its two vertices.

We call a line a {\em type $i,j$ line}
if one halfplane determined by that line has $i$ white vertices and
$j$ black vertices. In the other halfplane, there are $n-i-2$ white vertices and $n-j$
black vertices. Let $x_{i,j}$ be the
number of type $i,j$ lines.  Again, by symmetry, we can assume that
$0\leq j \leq  \lfloor \frac{n}{2} \rfloor = N+1$, and that $i<j$.
We will justify these assumptions later on.

Note that $y_i$ is related to $x_{i,j}$ by the following equation:
\begin{equation}\label{eqn_7}
y_i=x_{i,i+1}+\displaystyle\sum_{j=i+2}^{N+1} x_{i,j}
\end{equation}

Now, for a type $i,j$ line, there are  $i(n-j)+j(n-i-2)$ {\bf bw}-pairs of vertices
in opposite halfplanes of that edge. Summing this
quantity over all lines of a drawing, we obtain
$$M=\displaystyle\sum_{i=0}^{N}\displaystyle\sum_{j=i+1}^{N+1}
[i(n-j)+j(n-i-2)]x_{i,j}$$
{\bf bw}-pairs in opposite halfplanes.
As in the previous proof, we are looking for a drawing which minimizes $M$.

\medskip

We now justify our assumption that $0 \leq i  \leq N $ and $0  \leq j \leq N+1$
in a drawing which minimizes $M$. Assume that for a
given type $i,j$ line, the $i$ white vertices and the $j$ black
vertices are in different halfplanes. This yields $ij +
(n-j)(n-i-2)$ {\bf bw}-pairs.  However,
$$i(n-j)+j(n-i-2) \leq ij + (n-j)(n-i-2)$$
for $0 \leq i \leq N$ and $0 \leq j \leq N+1 $.
Therefore, the number of {\bf bw}-pairs over a drawing of
the graph is minimized when the $i$ and $j$ vertices are arranged so that
they lie in the same halfplane determined by the line.

Additionally, we justify our assumption that $i<j$, where $i$ and
$j$ respectively represent the number of white and black vertices in
one halfplane generated by two white points. Assume that in the
same halfplane, we had more white vertices than black vertices, i.e.
assume that $i$ is the number of black vertices and $j$ is the number of
white vertices, where $i<j$. Then we have $j(n-i) + i(n-j-2)$
{\bf bw}-pairs in opposite halfplanes. But, for $i<j$,
$$i(n-j)+j(n-i-2)< j(n-i) + i(n-j-2).$$
Therefore, the number of {\bf bw}-pairs over a drawing of the graph is minimized when a
halfplane determined by two white vertices contains more black
vertices (on the smaller side) than white vertices.

\medskip

In order to minimize $M$, we start by multiplying Equation (\ref{eqn_7}) by
$i(n-i-1)+ (i+1)(n-i-2)$, and subtracting it from $M$ for all values
of $i$, yielding:

\begin{eqnarray*}
& &M- \sum_{i=1}^{N} [i(n-i-1)+(i+1)(n-i-2)]y_i =\\
&=&\sum_{i=0}^{N}\sum_{j=i+1}^{N+1} \left[ i(n-j)+j(n-i-2) \right] x_{i,j}\\
& &\qquad -\sum_{i=1}^{N} [i(n-i-1)+(i+1)(n-i-2)]\left( \sum_{j=i+1}^{N+1} x_{i,j} \right)=\\
& =&\sum_{i=0}^{N}\sum_{j=i+1}^{N+1} \left[ i(n-j)+j(n-i-2) -i(n-i-1)-(i+1)(n-i-2)\right] x_{i,j}=\\
& =&\sum_{i=0}^{N}\sum_{j=i+1}^{N+1} \left[ i(i+1-j)+(j-i-1)(n-i-2) \right] x_{i,j} =\\
& =&\sum_{i=0}^{N}\sum_{j=i+1}^{N+1} (j-i-1)(n-2i-2) x_{i,j} = \\
&= &\sum_{i=0}^{N-1}\sum_{j=i+2}^{N+1} (j-(i+1))(n-2(i+1)) x_{i,j}
\end{eqnarray*}

Then:
\begin{eqnarray}\label{eqn_8}
M & = & \sum_{i=0}^{N} [i(n-i-1)+ (i+1)(n-i-2)]y_i +\\
\nonumber  & & + \sum_{i=0}^{N} \sum_{j=i+2}^{N+1} [(j-(i+1))(n-2(i+1))] x_{i,j}
\end{eqnarray}

As in the previous part, let $p_{s,t}$ be the number of white vertices having white endvertices
of type $s$ as the smallest type ($0 \leq s \leq N$) and  the index
$t$ is the number of distinct sorted sequences generated by all the vertices.  For example, in a
convex drawing of $K_{n,n}$, $p_{0,1}=n$, $p_{0,t}=0$ for $t \geq 2 $,
and $p_{s,t}=0$ for $s \geq 1 $.  Since $p_{s,t}$ counts white vertices and there are $n$ white vertices, we have:
$$n= \displaystyle\sum_{s=0}^{N} \displaystyle\sum_{t \geq 1} p_{s,t}.$$

Note that all the points on the convex hull have one distinct
sequence of endvertex types.

As in the previous part, denote by $z_{s,t,i}$ the number of appearances of $i$
in the $t$th sequence whose lowest term is $s$.  It follows that:
\begin{equation}\label{eqn_9}
y_i= \sum_{t \geq 1} \sum_{s=0}^{i} z_{s,t,i}p_{s,t},
\end{equation}
and for even $n$, we have
\begin{equation}\label{eqn_10}
y_N=\sum_{t \geq 1} \sum_{s=0}^{N} z_{s,t,N}p_{s,t}.
\end{equation}

Additionally, since there are $n$ white vertices, for a fixed $s$
and $t$, we have that:
\begin{equation}\label{eqn_11}
\displaystyle\sum_{i=s}^{N} z_{s,t,i}=n-1.
\end{equation}

Using the equality $\sum\limits_{t\geq 1} \sum\limits_{s=0}^{N} p_{s,t}=n$, Equations (\ref{eqn_9}) and (\ref{eqn_10}) become:
$$y_i= 2n + \sum_{t \geq 1}\left[ \sum_{s=0}^{i} (z_{s,t,i}-2)p_{s,t} -
2 \sum_{s=i+1}^{N}p_{s,t}\right],$$
and for even $n$, we have:
$$y_N=n+ \sum_{t \geq 1} \sum_{s=0}^{N} (z_{s,t,N}-1)p_{s,t}.$$

At this stage, one can easily see (similar to the proofs of Lemma \ref{lem_part1} above and
Lemma \ref{lem_part2} below) that if $s=0$,
the corresponding vertex is on the convex hull generated by the white points.
Hence, we have that $z_{0,t,i}=2$, and for even $n$, we have: $z_{0,t,N}=1$.
Therefore, since the first summand (for $s=0$) of the sum appearing in $y_i$ is $0$,
we can start the sum from $s=1$, without changing the value of $y_i$
(we actually could do it also in the first part, but it is unnecessary there). So, we have:
\begin{equation}\label{eqn_12}
y_i= 2n + \sum_{t \geq 1}\left[ \sum_{s=1}^{i} (z_{s,t,i}-2)p_{s,t} -
2 \sum_{s=i+1}^{N}p_{s,t}\right],
\end{equation}
and for even $n$, we have:
$$y_N=n+ \sum_{t \geq 1} \sum_{s=1}^{N} (z_{s,t,N}-1)p_{s,t}.$$

We proceed for odd $n$.  First note that:
$$i(n-i-1)+ (i+1)(n-i-2) =n(1+2i)- 2(i+1)^2.$$
Using Equation (\ref{eqn_12}), we can
rewrite the first part of the expression for $M$ (Equation (\ref{eqn_8}))  as
follows:
\begin{eqnarray*}
\sum_{i=0}^{N} \left[i(n-i-1)+ (i+1)(n-i-2)\right]y_i & =&\\
& &\hspace{-207pt}=\sum_{i=0}^{N} \left[n(1+2i) - 2(i+1)^2\right]y_i = \\
& &\hspace{-207pt}=2n \sum_{i=0}^{N} \left[n(1+2i) - 2(i+1)^2\right] +\\
& &\hspace{-200pt}+\displaystyle\sum_{t \geq 1}
\displaystyle\sum_{i=0}^{N} \left[n(1+2i) -
2(i+1)^2\right]\left[\displaystyle\sum_{s=1}^{i} (z_{s,t,i}-2)p_{s,t} -
2\displaystyle\sum_{s=i+1}^{N}p_{s,t}\right].
\end{eqnarray*}

Following a change in the indices of the sums, this can be rewritten as:
\begin{eqnarray*}
& & 2n \sum_{i=0}^{N} \left[n(1+2i) - 2(i+1)^2\right] +\\
& & \quad +\sum_{t \geq 1} \sum_{s=1}^{N}
p_{s,t}\left[\sum_{i=s}^{N} \left[n(1+2i) -
2(i+1)^2\right](z_{s,t,i}-2)\right.\\
& &\left.\quad\quad -2\sum_{i=0}^{s-1} \left[n(1+2i) -
2(i+1)^2\right]\right].
\end{eqnarray*}

This can again be rewritten as:
\begin{eqnarray*}
& &\hspace{-20pt} 2n \sum_{i=0}^{N} \left[n(1+2i) - 2(i+1)^2\right] +\\
& &\hspace{-20pt} \quad+\sum_{t \geq 1} \sum_{s=1}^{N} p_{s,t}\left[\left[n(1+2s) - 2(s+1)^2\right]\sum_{i=s}^{N} (z_{s,t,i}-2)+ \right.\\
& &\hspace{-20pt} \quad+ \sum_{i=s+1}^{N}\left(\left[n(1+2i) - 2(i+1)^2\right]-\left[n(1+2s) - 2(s+1)^2\right]\right)(z_{s,t,i}-2)\\
& &\hspace{-20pt} \quad \left.-2\sum_{i=0}^{s-1}\left[n(1+2i) - 2(i+1)^2\right]\right].
\end{eqnarray*}

Using Equation (\ref{eqn_11}), it follows that this is also equal to:
\begin{eqnarray*}
& &\hspace{-20pt} 2n \sum_{i=0}^{N}[n(1+2i) - 2(i+1)^2] + \left.\sum_{t \geq 1} \sum_{s=1}^{N} p_{s,t}\right[C(s,n)+ \\
& &\hspace{-20pt} \quad +\left.\sum_{i=s+1}^{N}\left(\left[n(1+2i) - 2(i+1)^2\right]-\left[n(1+2s) - 2(s+1)^2\right]\right)(z_{s,t,i}-2)\right],
\end{eqnarray*}
where
\begin{eqnarray*}
C(s,n)&=& \left[n(1+2s) - 2(s+1)^2\right]\left((n-1)- \displaystyle\sum_{i=s}^{N} 2\right) \\
& & -2\displaystyle\sum_{i=0}^{s-1} \left[n(1+2i) - 2(i+1)^2\right] =\\
&=& \left[n(1+2s) - 2(s+1)^2\right]((n-1)-2(N-s+1))\\
& &- 2\displaystyle\sum_{i=1}^{s-1} \left[n(1+2i) - 2(i+1)^2\right] - 2(n-2).
\end{eqnarray*}

We now show that $C(s,n)$ is non-negative for all $1\leq s \leq N$ and $n$.

\begin{lem}
For all $1\leq s \leq N$ and $n$, $C(s,n) \geq 0$.
\end{lem}

\begin{proof}
By the definition of $N=\lfloor \frac{n-1}{2} \rfloor$, we have
$$\left[n(1+2s) - 2(s+1)^2\right]((n-1)-2(N-s+1)) \geq \left[n(1+2s) - 2(s+1)^2\right](2s-1).$$
Then,
\begin{eqnarray*}
2\displaystyle\sum_{i=1}^{s-1} [n(1+2i) - 2(i+1)^2] &\leq& 2\displaystyle\sum_{i=1}^{s-1}[n(1+2(s-1)) - 2((s-1)+1)^2]=\\
&=& 2\displaystyle\sum_{i=1}^{s-1}[n(2s-1) - 2s^2]= \\
&=& 2(s-1)[n(2s-1) - 2s^2]
\end{eqnarray*}

Therefore,
\begin{eqnarray*}
C(s,n) & \geq &  (2s-1)\left(n(1+2s) - 2(s+1)^2\right) \\
& &- 2(s-1)\left(n(2s-1) - 2s^2\right) -2(n-2)=\\
&=&-2s^2 +6ns-5n+6 =\\
&=& (-2s^2+ns)+(5ns-5n)+6 \geq 0,
\end{eqnarray*}
for all $1 \leq s \leq N$.
\end{proof}

We now show that $z_{s,t,i} -2 \geq 0$: Since the term $z_{s,t,i} -1$ must
be carried throughout this summation, the expression for even $n$ is
also minimized for $p_{s,t}=0$ for all $s \geq 1$, provided $z_{s,t,N} -1 \geq 0$.

\begin{lem}\label{lem_part2}
$z_{s,t,i} \geq 2$ for all $s,t,i$ such that $i>s$, and $z_{s,t,N} \geq 1$ for even $n$.
\end{lem}

\begin{proof}
For a given white vertex, we start by proving that there is at least one line with an endvertex of type $\frac{n-2}{2}$ for even $n$ and at least two lines with endvertices of type $\frac{n-3}{2}$ for odd $n$.  This statement can be proved by induction on $n$.  This statement is obvious for $n=2$ and $n=3$, so we start with the inductive step.  Also, note that in traversing the $n$ lines  incident to a given vertex in a clockwise or counterclockwise manner in moving from line to line, the number of white vertices in the clockwise following halfplane may be changed by at most $1$. This fact will be used numerous times throughout the proof (see Example \ref{example1}).

\medskip

\noindent
\textbf{Case I:} Passing from even $n$ to $n+1$. \\
We consider the line with an endvertex of type $\frac{n-2}{2}$ in the
configuration of $n$ white vertices and $n$ black vertices. When the
$(n+1)$st pair of vertices is added, the original endvertex will be the first
vertex of type $\frac{(n+1)-3}{2}$. If the $(n+1)$st white vertex is
added in this line's clockwise following halfplane, then an endvertex
of an immediately following line has type
$\frac{(n+1)-3}{2}$.  Thus, the endvertex of this line will be the second endvertex of
type $\frac{n-2}{2}$.

\medskip

\noindent
\textbf{Case II:} Passing from odd $n$ to $n+1$. \\
Consider a line with an endvertex of type $\frac{n-3}{2}$, which has
$\frac{n-3}{2}$ white vertices in one of its halfplanes and
$\frac{n-1}{2}$ in the other.  If the $(n+1)$st white vertex is added
in the halfplane with $\frac{n-3}{2}$ white vertices, then the
endvertex of the considered line is of type $\frac{(n+1)-2}{2}$.  If
the $(n+1)$st white vertex is added in the halfplane with
$\frac{n-1}{2}$ white vertices, then there are $\frac{n+1}{2}$ white
vertices in this halfplane and $\frac{(n+1)-4}{2}$ white vertices in
the clockwise following halfplane of this line.  Since
the number of white vertices in the clockwise following halfplane
can be changed by at most $1$ while moving from line to line,
we find that traversing the graph
from the line with $\frac{n+1}{2}$ white vertices in the clockwise
following halfplane to the line with $\frac{(n+1)-4}{2}$, there
must occur a line with $\frac{(n+1)-2}{2}$ white
vertices in the clockwise following halfplane.  Thus, the endvertex
of this line is of type $\frac{(n+1)-2}{2}$.

\medskip

Using this result and the fact that while moving
from line to adjacent line, the number of white vertices
in the clockwise following halfplane may be changed by at most $1$, we
can prove that there are two lines with endvertices of each type
from the type $s+1$ to the maximal type $N$.

\medskip

We split the proof according to the parity of $n$.

\begin{itemize}
\item For even $n$, we have one line with an endvertex of maximal type $N=\frac{n-2}{2}$.
Traversing the $n-1$ lines starting and ending with the line with
an endvertex of type $N$ from line to line, we must go down to
a line with $s$ white vertices in the clockwise
following halfplane, and then back up to one linewith $N$.  Thus, we
find that there are at least two lines with endvertices of
each type from $s+1$ to $N$.

\item For odd $n$, we have two lines of
maximal type $N=\frac{n-3}{2}$.  Traversing the $n-1$ lines from one line
with endvertex of type $N$ to the other, we must go down to a line
with $s$ white vertices in the clockwise following
halfplane and back up to one line with $N$.  Thus, there are at least two
lines with endvertices of each type from $s+1$ to $N$.
\end{itemize}

\medskip

Hence, it follows that $z_{s,t,i} \geq 2$ for all $s,t,i$ such that $i>s$, and $z_{s,t,N} \geq 1$ for even $n$, as needed.
\end{proof}

Additionally, by Observation \ref{obs1},
\begin{eqnarray*}
\sum_{i=s+1}^{N}([n(1+2i) - 2(i+1)^2]-[n(1+2s) -
2(s+1)^2]) &=&\\
=\sum_{i=s+1}^{N} 2[i(n-i-2)-s(n-s-2)] & \geq & 0
\end{eqnarray*}
for $i,s \leq N=\frac{n-1}{2}.$

Going back to the final expression for Equation (\ref{eqn_8}), we have:

\begin{eqnarray*}
\hspace{-10pt}M&=&2n \sum_{i=0}^{N} \left[n(1+2i) - 2(i+1)^2\right] + \left. \sum_{t \geq 1} \sum_{s=1}^{N} p_{s,t} \right[C(s,n)+\\
\hspace{-10pt}& & \left.+\sum_{i=s+1}^{N}\left(\left[n(1+2i) - 2(i+1)^2\right]-\left[n(1+2s) - 2(s+1)^2\right]\right)(z_{s,t,i}-2)\right] \\
\hspace{-10pt}& & +\sum_{i=0}^{N} \sum_{j=i+2}^{N+1} [(j-(i+1))(n-2(i+1))] x_{i,j}
\end{eqnarray*}

Since $C(s,n)$, $(z_{s,t,i} -2)$ and $(j-(i+1))(n-2(i+1))$ are
non-negative, we find that this expression is
minimized when $p_{s,t}=0$ for $s \geq 1$ and $x_{i,j}=0$ for $i<j$.
Evaluating the sum for these conditions, we
have:
$$M=2n \displaystyle\sum_{i=0}^{N} \left[n(1+2i) - 2(i+1)^2\right] = 2n{n \choose 3}.$$
Since each line determined by two white points was
counted twice, once for each endvertex, we have $n{n \choose 3}$ {\bf bw}-pairs in opposite halfplanes, as claimed.

As in the previous part, if we perform similar computations for even $n$, we get the same result.
\end{proof}

\subsection{Final step of the proof}\label{third}

Here, we finish the proof of Theorem \ref{main_thm}.

\begin{proof}[Proof of Theorem \ref{main_thm}] For a given rectilinear drawing $D$ of $K_{n,n}$, any
Orchard crossing is determined by a {\bf bw}-pair, where
one point is in one halfplane of line $\ell_k$ and the other point is
in the other halfplane of $\ell_k$, where the type of the line $\ell_k$ is one of the following three types:
\begin{itemize}
\item[(a)] a line determined by a {\bf bw}-pair.
\item[(b)] a line determined by two white points.
\item[(c)] a line determined by two black points.
\end{itemize}

Let $A,B,C$ be the numbers of {\bf bw}-pairs determined by lines of types (a),(b),(c), respectively.
Then, $n(D) = A+B+C$. By Propositions \ref{A} and \ref{B}, we have $A \geq 2n{n \choose 3}$ and $B,C \geq n{n \choose 3}$.
Hence:
$$n(D) \geq 2n{n \choose 3} + n{n \choose 3} + n{n \choose 3}= 4n{n \choose 3}$$

On the other hand, by Lemma \ref{num_comp_bipar}, we have a drawing $D$ of $K_{n,n}$ which satisfies $n(D)=4n {n \choose 3}$.
So finally we have ${\rm OCN}(K_{n,n})=4n {n \choose 3}$ as claimed.
\end{proof}

\end{document}